\documentclass[11pt]{article}

\usepackage{amsmath}
\usepackage{amssymb}
\usepackage{bbm}
\usepackage{amsfonts}
\usepackage{amsthm}
\usepackage{bbm}
\usepackage{parskip}
\usepackage{tikz}
\usepackage{fullpage}

\newtheorem{thm}{Theorem}[section]
\newtheorem{prop}[thm]{Proposition}

\newtheorem{lem}[thm]{Lemma}

\title{Uniform Oscillatory Integral estimates for Convex Phases via Sublevel Set estimates}
\author{John Green}
\date{}

\begin{document}
\maketitle

\begin{abstract}
We examine the relation between oscillatory integral estimates and sublevel set estimates associated to convex functions. Whilst the former implies the latter in many cases, the reverse requires additional assumptions. Under finite (line) type assumptions, Bruna, Nagel \& Wainger \cite{bruna1988convex} were able to demonstrate a very precise control of oscillatory integrals with convex phases via their sublevel sets. Without the finite type assumption, certain erratic behaviour can force this precise control to fail (Bak, McMichael, Vance \& Wainger \cite{bak1989fourier}). We establish the same precise control under an alternative qualitative geometric assumption.
\end{abstract}

\begin{footnotesize}
\end{footnotesize}
\section{Introduction}
We consider oscillatory integrals of the form
$$I(\lambda)=\int_\Omega \eta(x)e^{i\lambda f(x)}\,dx$$
where $\Omega$ is an open, bounded convex subset of $\mathbb{R}^n$, $f$ is a $C^2$ convex function on that set (the ``phase") with a unique minimum in $\Omega$, $\lambda$ is a real ``frequency parameter", and $\eta$ is an ``amplitude", which we shall generally take to be $1$. Note that by convexity, the unique minimum is the unique critical point. We will be interested in bounds for $|I(\lambda)|$, and thus we may assume by translation of the phase that the value of $f$ at the critical point is $0$.

It is well known that the behaviour of such oscillatory integrals is largely governed by the behaviour of the phase near critical points, and in particular, how quickly one moves away from these critical points, as determined by higher order derivatives of the phase (see, for instance, Stein \cite{stein1993harmonic}).

This leads one to consider the measures of the sets of points near those critical points, that is, the sublevel sets
$$S(\varepsilon)=\{x\in\Omega:|f(x)|\leq\varepsilon\}.$$
In many situations, it is known that oscillatory integral estimates imply sublevel set estimates in a way that exhibits a certain kind of uniformity. More precisely, we have

\vspace*{0.3cm}

\begin{prop}\label{oscimpsub}
Let $f:\Omega\rightarrow\mathbb{R}$ be a measurable function such that for all real non-zero $\lambda$, we have
$$\left|\int_\Omega e^{i\lambda f(x)}\,dx\right|\leq A|\lambda|^{-\delta}$$
where $0<\delta<1$. Then for each $c\in\mathbb{R}$, we have
$$|\{x\in\Omega:|f(x)-c|\leq\varepsilon\}|\leq C_\delta A\varepsilon^\delta$$
where $C_\delta$ depends only on $\delta$.
\end{prop}

See, for instance, Carbery-Christ-Wright \cite{carbery1999multidimensional}.

We will be concerned with going in the opposite direction, more precisely, given a suitable increasing function $t$ with $t(0)=0$, can we show
$$|S(\varepsilon)|\leq K_1t(\varepsilon)\Rightarrow|I(\lambda)|\leq K_2t(1/|\lambda|)$$
with $K_2$ independent of $f$ and depending boundedly on $K_1$ (and $\eta$, $n$ and $\Omega$) \footnote{ In Proposition \ref{oscimpsub}, we considered sublevel sets at each height $c$; here we focus only on the height of the critical value. The main reason for this is that the principal contribution to the oscillatory integral comes from the critical points, but it should also be noted that Proposition \ref{oscimpsub} fails for $\delta\geq 1$, and in dimensions greater than one, sublevel set estimates with $\delta>1$ can only hold at critical values.}?

Bruna, Nagel \& Wainger \cite{bruna1988convex} obtained a result that gives very precise control of the oscillatory integral by the sublevel sets in the finite (line) type setting - roughly speaking, $f$ is of finite line type $m$ if when restricted to each line its derivatives do not vanish to order $m$. In this setting they showed
$$|I(\lambda)|\leq C|S(1/|\lambda|)|$$
for fixed $\eta\in C_c^\infty(\Omega)$, where the constant $C$ depends on $m$ and quantitative properties of $f$ in a complicated way, but these dependencies are such that the constant remains bounded for certain families of phase functions $f$. In particular, they are sufficient for deriving certain important consequences, such as Fourier transform estimates for surface measures of convex hypersurfaces of finite line type. See also Cowling, Disney, Mauceri \& M\"uller \cite{cowling1990damping}.

In addition to this result, Bruna, Nagel \& Wainger proved a similar result without the finite type assumption by a much simpler argument in the case where $\Omega$ is an interval in $\mathbb{R}$. This argument also gives a much simpler constant. One might hope to extend this to higher dimensions, but this turns out to be impossible even when $\Omega\subseteq\mathbb{R}^2$ and $f$ is radial, as was shown by Bak, McMichael, Vance \& Wainger \cite{bak1989fourier}.

At the heart of the issue is this: An estimate on $|S(1/|\lambda|)|$ says nothing about how $f$ behaves outside of that set; even with the convexity assumption, the growth of $|S(\varepsilon)|$ can be erratic - for instance, one can imagine in the radial setting that $f$ grows slowly and linearly, then increases rapidly for a short time, then linearly again. Naturally, since oscillatory integrals depend on how $f$ behaves everywhere, we shouldn't expect Bruna-Nagel-Wainger type estimates to hold without something that implicitly controls all the other sublevel sets, and convexity alone proves insufficient.

The finite type assumption provides this extra information; however, if we are not interested in the precise control of the Bruna-Nagel-Wainger estimate and instead want to skip straight to its consequence, that estimates $|S(\varepsilon)|\leq K_1t(\varepsilon)$ imply estimates $|I(\lambda)|\leq K_2t(|\lambda|^{-1})$, we note that assuming the former already contains some information on how $f$ behaves everywhere, in the form of an estimate for \textit{each} sublevel set. Of course, this does not allow for $t$ to be an arbitrary increasing function, for then setting $t(x)=|S(x)|$, the Bak-McMichael-Vance-Wainger counterexample would still apply. The following condition on $|S(x)|$ will prove to be sufficient for our purposes:

\textbf{Geometric assumption.} Let $l$ be such that $S(l)$ is compact. We assume that on $[0,l]$, $|S(y)|$ is a concave function of $y$. Equivalently, for each $h\in(0,l)$, $|\{x\in\Omega:y\leq f(x)\leq y+h\}|$ is a decreasing function of $y$ on $[0,l-h]$.

The equivalence of these assumptions can be seen through the coarea formula, we give the details later. We will also note that for convex functions with some $S(l)$ compactly contained in $\Omega$, $|S(y)|^{1/n}$ is necessarily  concave, so that the assumption is always true for $n=1$, matching the observation of Bruna, Nagel \& Wainger.

\textit{Remark.} One can check that this is true for many convex functions that vanish to infinite order at their minimum. However, there are also functions of finite type for which this is false, so this does not generalise the Bruna-Nagel-Wainger estimate.

Additionally, the alternative formulation of the assumption establishes a connection with recent work of Basu, Guo, Zhang \& Zorin-Kranich \cite{basu2021stationary}, wherein they bound oscillatory integrals by sublevel sets under the assumption that for each fixed width the measure of the sublevel sets as a function of height changes monotonicity boundedly many times, and establish bounds on the number of changes of monotonicity for semialgebraic functions of bounded complexity. Their method possesses similarities with ours and can be adapted without much difficulty to give weaker versions of our results in the present setting.

\textbf{Notations. }The Euclidean norm of $x$ will be denoted by $|x|$, and the Lebesgue and Hausdorff measures of sets $A$ will be denoted $|A|$, exactly which measure is intended will be clear from the context. The gradient vector is a row vector denoted by $\nabla f$, and the Hessian will be denoted $\mathcal{H}f$. The Jacobian matrix of a diffeomorphism $\psi$ will be denoted $d\psi$. The standard basis of $\mathbb{R}^n$ will be denoted $e_i$, and the volume of the unit ball in $\mathbb{R}^n$ denoted by $\omega_n$. The tangent space to a manifold $M$ at point $x$ is denoted $T_xM$.

\section{Main results}
Throughout the rest of the paper we will use the following set-up: $\Omega$ shall be an open, bounded convex subset of $\mathbb{R}^n$, $f$ a $C^2$ convex function on that set with a unique minimum $0$ (hence unique critical point with critical value $0$) in $\Omega$, and the geometric assumption is satisfied, that is, there is $l$ such that $S(l)$ is compact, and on $[0,l]$, $|S(y)|$ is a concave function of $y$. Let $\eta$ have the form $a(f(x))$ for some $C^1$ function on $[0,l]$, that is, $\eta$ is constant on level sets of $f$.

Let $I(\lambda)$ denote the oscillatory integral associated to $f$,
$$I(\lambda)=\int_{S(l)} \eta(x)e^{i\lambda f(x)}\,dx.$$

\vspace*{0.3cm}

\begin{thm}
Subject to the above assumptions, we have
$$|I(\lambda)|\leq 5n(\|a\|_{L^\infty[0,L]}+\|a'\|_{L^1[0,L])})|S(\omega_n^{-1/n}/|\lambda|)|.$$
\end{thm}

This result as stated is intended to express concisely the kind of consequence that follows from our method. This method developed from an investigation into the role of geometric structure in oscillatory integral estimates, started in an earlier work of the author \cite{green2021algebraic}, and was not developed with any particular application in mind. As such, the statement above is not formulated for the greatest applicability, rather, we stress that our method should be suitably adapted once a particular application is in mind. We provide here some remarks to facilitate this, but stress that these comments are best understood after seeing the proofs in the following sections.

\begin{enumerate}
\item The result remains true if all the measures involved in the assumptions and in the conclusions are replaced with weighted measures with a sufficiently smooth positive weight. The arguments proceed identically. In particular, if one is interested in surface measure on the graph of $f$, one may use the weight $\sqrt{1+|\nabla f|^2}$.
\item If $a$ is monotone, then $\|a'\|_1\leq\|a\|_\infty$.
\item If one wishes to consider an amplitude function $\eta$ not constant on level sets of $f$, one may define $$a(s)=\left(\int_{f^{-1}(s)}\frac{\eta(x)}{|\nabla f(x)|}\,d\mathcal{H}^{n-1}(x)\right)/\left(\int_{f^{-1}(s)}\frac{1}{|\nabla f(x)|}\,d\mathcal{H}^{n-1}(x)\right).$$
Following arguments in section \ref{prelim}, we can show that this is $C^1$ as a function of $s$, and we can see $a(s)$ as a weighted average of $\eta$, from which it is clear that $\|a\|_\infty\leq\|\eta\|_\infty$. Finally, by the coarea formula, we have
\begin{align*}
\int_{S(l)} a(f(x))e^{i\lambda f(x)}\,dx&=\int_0^le^{i\lambda s}\int_{f^{-1}(s)}\frac{a(s)}{|\nabla f(x)|}\,d\mathcal{H}^{n-1}(x)\,ds\\
&=\int_0^le^{i\lambda s}\int_{f^{-1}(s)}\frac{\eta(x)}{|\nabla f(x)|}\,d\mathcal{H}^{n-1}(x)\,ds\\
&=\int_{S(l)} \eta(x)e^{i\lambda f(x)}\,dx
\end{align*}
Then if $a$ is monotone, by the second remark we can estimate by $C_n\|\eta\|_\infty|S(\omega_n^{-1/n}/|\lambda|)|$.
\item Furthermore, one can use the first and third remarks together, and factor an amplitude function into a function constant on level sets and a weight with respect to which the sublevel sets are concave.
\end{enumerate}

Finally, we note here an equivalent condition to the geometric assumption, which we shall prove during section \ref{prelim}.

\vspace*{0.3cm}

\begin{lem}\label{derivcond}
In the setting described, $|S(y)|$ is a concave function of $y$ on $[0,l]$ if and only if
$$\int_{f^{-1}(s)}\frac{|\nabla f|^2\Delta f- 2(\nabla f) \mathcal{H}f(\nabla f)^{T}}{|\nabla f|^5}\,d\mathcal{H}^{n-1}(x)\leq 0$$
for each $s\in(0,l)$. Moreover, in the radial setting $f(x)=F(|x|)$, this condition is equivalent to $(n-1)F'(s)\leq sF''(s)$.
\end{lem}

\section{Preliminaries}\label{prelim}
We begin by recalling the coarea formula for $f$ (see, for instance, Federer \cite{federer2014geometric}). For a measurable function $k(x)$ on $\Omega$, we have
$$\int_\Omega k(x)|\nabla f(x)|\,dx=\int_\mathbb{R} \int_{f^{-1}(s)} k(x)\,d\mathcal{H}^{n-1}(x)\,ds$$
whenever both sides of the equation exist, where $\mathcal{H}^{n-1}$ denotes the $(n-1)$-dimensional Hausdorff measure. Later, we shall apply this to the oscillatory integral to obtain a one-dimensional oscillatory integral, this technique has appeared previously in the literature, see Cowling, Disney, Mauceri \& M\"uller \cite{cowling1990damping} and the references therein, in particular, Varchenko \cite{varchenko1976newton}. However, we first apply the coarea formula with $k(x)=\chi_{\{a\leq f(x)\leq b\}}(x)|\nabla f(x)|^{-1}$, where $\chi$ denotes an indicator function, yielding
$$|\{x\in\Omega:a\leq f(x)\leq b\}|=\int_a^b\int_{f^{-1}(s)}\frac{1}{|\nabla f(x)|}\,d\mathcal{H}^{n-1}(x)\,ds.$$
In what follows we shall denote the integrand of the right hand side by $J(s)$ for convenience. Later we shall show that $J$ is a $C^1$ function of $s$. Assuming this momentarily, we see that the geometric assumption given in the introduction is in fact equivalent to $J'(s)\leq 0$. Indeed, setting $a=0$, $b=y$ in the above, we see that $|S(y)|$ is obtained from integrating $J(s)$ from $0$ to $y$, in particular, the derivative of $|S(y)|$ is $J(y)$. Since concavity of a function is equivalent to non-positivity of the second derivative, the equivalence follows.

Likewise, we also note that $J'(y)\leq 0$ is equivalent to $J(y)$ being non-increasing, this is equivalent to saying $J(y+h)-J(y)\leq 0$ for each $h$. Setting $a=y$, $b=y+h$, we see this quantity is the derivative of $|\{x\in\Omega:y\leq f(x)\leq y+h\}|$, hence $J'(y)\leq 0$ is equivalent to this being a decreasing function of $y$ for each $h$, as outlined in the introduction.

We will use non-decreasing Schwarz symmetrisation (see Kawohl \cite{kawohl2006rearrangements}), replacing $f$ with a radial non-decreasing rearrangement $\tilde{f}$. Define $g:[0,l]\rightarrow [0,(\omega_n^{-1}S(l))^{1/n}]$ by
$$g(y):=(\omega_n^{-1}|\{x\in\Omega:f(x)\leq y\}|)^{1/n},$$
where $\omega_n$ is the volume of the unit ball in $\mathbb{R}^n$. One can see that $g$ is $C^2$ on $(0,l)$ and strictly increasing on $[0,l]$ (we have seen above that $|\{x\in\Omega:f(x)\leq y\}|$ is $C^2$ with strictly positive derivative $J(y)$, from this the conclusion follows), hence $g^{-1}$ exists and is strictly increasing. We define $\tilde{f}(x)$ on the ball $B$ of radius $(\omega_n^{-1}S(l))^{1/n}$ by $\tilde{f}(x)=g^{-1}(|x|)$. Then
\begin{align*}
|\{x\in B:\tilde{f}(x)\leq y\}|&=|\{x\in B:g^{-1}(|x|)\leq y\}|=|\{x\in B:|x|\leq g(y)\}|=\omega_ng(y)^n\\
&=|\{x\in\Omega:f(x)\leq y\}|.
\end{align*}
Note that it follows from $g$ being $C^2$ that $\tilde{f}$ is $C^2$ away from the origin, and we can in fact apply the coarea formula in the same way. Since the function $|S(y)|$ associated to $\tilde{f}$ is identical to the one associated to $f$, so is its derivative, hence we can use the alternative expression
$$J(s)=\int_{\tilde{f}^{-1}(s)}\frac{1}{|\nabla \tilde{f}(x)|}\,d\mathcal{H}^{n-1}(x).$$

Now, it is in fact true that $g$ is concave - this can be shown using the Brunn-Minkowski inequality. It follows that $g^{-1}$ is convex, see Theorem $4.1$ of \cite{berman2014symmetrization}. Hence $\tilde{f}$ is convex on lines through the origin.

It remains to clarify that $J(s)$ is indeed $C^1$. In fact, we shall calculate the derivative in order to obtain Lemma \ref{derivcond}. First, note that since $0$ is the only critical value of $f$, the level sets $f^{-1}(s)$ are $C^2$ manifolds for $s\in(0,l)$. For fixed $s\in(0,l)$, let $\chi$ be a smooth bump function with compact support in $(0,l)$ and equal to $1$ in a neighbourhood of $s$. Let $X$ be a vector field on $\Omega$ defined by $X(x)=\chi(f(x))\nabla f(x)/|\nabla f(x)|^2$, this is $C^1$ since $f$ is $C^2$ with $\nabla f(x)\neq 0$ on the support of $\chi(f(x))$. We will be interested in the flow of this vector field within a small neighbourhood (see, for instance, Lee \cite{lee2013smooth}).

The vector field $X$ defines, for small $r$ in a neighbourhood of $0$, a one-parameter family of diffeomorphisms $\psi_r(x)=\psi(x,r)$ by the equation
$$\frac{\partial}{\partial r}\psi(x,r)=X(\psi(x,r))$$
in a neighbourhood of $f^{-1}(s)$. Note $\psi_r^{-1}=\psi_{-r}$, and $\psi_0$ is simply the identity. We will only be interested in these diffeomorphisms in a neighbourhood $U$ of $f^{-1}(s)$ small enough that $\chi(f(x))=1$ for all $x$ in $\psi_r(U)$.

It is easy to see that $\psi_h$ maps $f^{-1}(s)$ to $f^{-1}(s+h)$, since
$$\frac{\partial}{\partial r}f(\psi(x,r))=\nabla f(\psi(x,r))\cdot X(\psi(x,r))=1$$
by definition of $X$, and it follows that $\psi_h$ restricts to a diffeomorphism of $f^{-1}(s)$ to $f^{-1}(s+h)$. Using change of variables for functions between manifolds (for instance, Spivak \cite{spivak1975comprehensive}) it follows that
$$J(s+h)=\int_{f^{-1}(s)}\det\left(d\psi_h|_{T_xf^{-1}(s)\rightarrow T_{\psi_h(x)}f^{-1}(s+h)}\right)|\nabla f(\psi_h(x))|^{-1}\,d\mathcal{H}^{n-1}(x).$$
We shall express this determinant in terms of the determinant of $d\psi_h$. To do this, note that the orthogonal component to each of these tangent spaces is given by the gradient, so we can choose orthogonal rotations $A$ and $B$ so that $A(e_n)=\nabla f(x)/|\nabla f(x)|$ and $B(\nabla f(\psi_h(x))/|\nabla f(\psi_h(x))|)=e_n$. Then $Bd\psi_h(x)A$ has the same determinant as $d\psi_h(x)$, and is such that the $n^{\text{th}}$ row is $0$ besides the $n^{\text{th}}$ column, call this entry $R$, and the first $(n-1)\times(n-1)$ rows and columns represent $d\psi_h|_{T_xf^{-1}(s)\rightarrow T_{\psi_h(x)}f^{-1}(s+h)}$ in an orthonormal basis. Hence
$$\det\left(d\psi_h|_{T_xf^{-1}(s)\rightarrow T_{\psi_h(x)}f^{-1}(s+h)}\right)=\frac{\det(d\psi_h(x))}{R}.$$
It remains to calculate $R$, that is the component of $d\psi_h(x)(\nabla f(x)/|\nabla f(x)|)$ in the direction of $\nabla f(\psi_h(x))/|\nabla f(\psi_h(x))|$. Now, $\gamma(t)=\psi(x,t)$ is a short curve with $\gamma(0)=x$ and $\gamma'(0)=\nabla f(x)/|\nabla f(x)|^2$, so we have
\begin{align*}
d\psi_h(x)(\nabla f(x)/|\nabla f(x)|)&=|\nabla f(x)|d\psi_h(x)(\nabla f(x)/|\nabla f(x)|^2)\\
&=|\nabla f(x)|(\psi_h\circ\gamma)'(0)=|\nabla f(x)|(\psi(x,t+h))'(0)\\
&=|\nabla f(x)|\frac{\nabla f(\psi_h(x))}{|\nabla f(\psi_h(x))|^2}.
\end{align*}
It follows from this that $R=|\nabla f(x)|/|\nabla f(\psi_h(x))|$, and thus
$$J(s+h)=\int_{f^{-1}(s)}\frac{\det(d\psi_h(x))}{|\nabla f(x)|}\,d\mathcal{H}^{n-1}(x).$$
One may now differentiate under the integral sign to obtain
$$J'(s)=\int_{f^{-1}(s)}\frac{[\text{div}(X)](x)}{|\nabla f(x)|}\,d\mathcal{H}^{n-1}(x)$$
where we have used the non-linear version of the Liouville-Ostrogradski formula (see, for instance, Teschl \cite{teschl2012ordinary}) to obtain $(d/dh)|_{h=0}\det(d\psi_h(x))=[\text{div}(X)](x)$. The divergence of $f$ may be calculated as
$$[\text{div}(X)](x)=\frac{|\nabla f|^2\Delta f-2(\nabla f)\mathcal{H}f(\nabla f)^{T}}{|\nabla f|^4}$$
from which it is clear that $J'(s)$ is continuous (by using the same change of variables as before to express each integral as an integral over the same level set), and the first part of Lemma \ref{derivcond} follows.

The reduction to the radial case can be calculated directly from this, however, it is easier in that case to simply note that the condition that $|S(y)|$ is concave is equivalent to, in the terminology of the lemma, $(F^{-1}(s))^n$ being a concave function. Since concavity is equivalent to non-positivity of the second derivative, one sees that this is in turn equivalent to
$$n(n-1)(F^{-1}(s))^{n-2}((F^{-1})'(s))^2+n(F^{-1}(s))^{n-1}(F^{-1})''(s)\leq 0$$
holding for each $s$. Using $(F^{-1})'(s)=1/F'(F^{-1}(s))$ and $(F^{-1})''(s)=-F''(F^{-1}(s))/((F^{-1})'(s))^3$, one rearranges to obtain the desired result.

\section{Proof of the main theorem}
The proof goes as follows. We split the integral over a certain sublevel set, at height $\varepsilon_0$, say, and its complement. To the former we apply trivial bounds and estimate by the measure of the sublevel set, for the latter we integrate by parts as standard in many oscillatory integral arguments, and estimate by $J(\varepsilon_0)$. To estimate $J(\varepsilon_0)$, we use its equivalent form arising from the Schwarz symmetrisation and estimate by the measure of the level set times $1/|\nabla\tilde{f}|$. Since the level sets of $\tilde{f}$ are spheres, we can write this as a function of the measure of its interior, which is $|S(\varepsilon_0)|$. To bound the gradient of $\tilde{f}$, we exploit the radial nature to relate the diameter of a sublevel set to its measure, and obtain sublevel set inclusions for $|\nabla \tilde{f}|$ in terms of sublevel sets of $\tilde{f}$, generalising an observation in one dimension used in an earlier work of the author \cite{green2021algebraic}.

We proceed via the coarea formula. For $\varepsilon_0\in(0,l)$ we have
$$\int_\Omega \eta(x)e^{i\lambda f(x)}\,dx=\int_{S_{\varepsilon_0}}\eta(x)e^{i\lambda f(x)}\,dx+\int_{\varepsilon_0}^le^{i\lambda s}a(s)\int_{f^{-1}(s)}\frac{1}{|\nabla f(x)|}\,d\mathcal{H}^{n-1}(x)\,ds.$$
The first term can be estimated by $\|a\|_\infty |S(\varepsilon_0)|$. To estimate the latter, write $e^{i\lambda s}=(i\lambda)^{-1}(d/ds)e^{i\lambda s}$ and integrate by parts. We have
\begin{align*}
\int_{\varepsilon_0}^le^{i\lambda s}a(s)J(s)\,ds=(i\lambda)^{-1}\left[e^{i\lambda l}J(l)a(l)-e^{i\lambda \varepsilon_0}J(\varepsilon_0)a(\varepsilon_0)-\int_{\varepsilon_0}^le^{i\lambda s}\left(\frac{d}{ds}\right)a(s)J(s)\,ds\right]
\end{align*}
We use the triangle inequality and estimate each term. The first two terms are bounded by $\|a\|_\infty J(\varepsilon_0)$, since $J$ is non-increasing. For the integral, we have
\begin{align*}
\left|\int_{\varepsilon_0}^le^{i\lambda s}\left(\frac{d}{ds}\right)a(s)J(s)\,ds\right|&\leq\int_{\varepsilon_0}^l|a'(s)J(s)|\,ds+\int_{\varepsilon_0}^l|a(s)J'(s)|\,ds\\
&\leq J(\varepsilon_0)\|a'\|_1+\|a\|_\infty\int_{\varepsilon_0}^l|J'(s)|\,ds
\end{align*}
and use that $J'$ is single signed to take the absolute value signs out the integral and estimate the integral by its endpoints, giving a bound of $2J(\varepsilon_0)$. Altogether, we get
$$\left|\int_{\varepsilon_0}^le^{i\lambda s}a(s)J(s)\,ds\right|\leq 4(\|a\|_\infty+\|a'\|_1)|\lambda|^{-1}J(\varepsilon_0).$$

Now, denote $t(x)=|S(x)|$. By definition we have $|\{x\in B:\tilde{f}(x)\leq r\}|=t(r)$, and we know that the radial derivative $\tilde{f}_r$ of $\tilde{f}$ is increasing, since, as noted earlier, $\tilde{f}$ is convex along lines through the origin.

Define $A_\alpha=\{x\in B\setminus\{0\}:|\tilde{f}_r|=|\nabla\tilde{f}(x)|\leq \alpha\}$. Since $\tilde{f}(x)$ is radial, so is $A_\alpha$. Moreover, since $\tilde{f}_r$ is increasing, there is an $L$ so that $|\nabla f(x)|$ is at most $\alpha$ for $|x|\leq L$ and $x$ is not in $A_\alpha$ otherwise. Thus its measure is equal to $\omega_nL^n$. Now, if for some $x\in A_\alpha$ we have $\tilde{f}(x)=c>0$, then by the mean value theorem we have $c/L\leq \tilde{f}(x)/|x|=|\tilde{f}(x)-\tilde{f}(0)|/|x-0|=\tilde{f}_r(s)$ for some $s<|x|\leq L$. But $\tilde{f}_r(s)\leq\alpha$, hence $c\leq L\alpha$. Thus $A_\alpha$ is contained in the set $\{x\in B:\tilde{f}(x)\leq L\alpha\}$, and so by the sublevel set estimate, we have $\omega_nL^n=|A_\alpha|\leq t(L\alpha)$.

Now, since $t$ is strictly increasing and thus has strictly increasing inverse, we have $t^{-1}(\omega_nL^n)/L\leq\alpha$. Define a new function $T$ by $T(x)=t^{-1}(\omega_nx^n)/x$. We would like to check that this is strictly increasing. We calculate the derivative
$$T'(x)=(n\omega_n)(t^{-1})'(\omega_nx^n)x^{n-2}-\frac{t^{-1}(\omega_nx^n)}{x^2}$$
which is positive for each $x>0$ if and only if $(n\omega_n)(t^{-1})'(\omega_nx^n)x^n>t^{-1}(\omega_nx^n)$ for $x>0$. Writing $x=((\omega_n)^{-1}t(y))^{1/n}$, and using $(t^{-1})'(z)=1/t'(t^{-1}(z))$, this inequality becomes $nt(y)/t'(y)>y$.

Now, in the notation of the preliminaries, note that $t(y)=\omega_n(g(y))^n$ for the strictly increasing concave function $g$. It follows that $nt(y)/t'(y)>y$ is equivalent to $g(y)/g'(y)>y$, which is equivalent to the strictly increasing convex function $g^{-1}$ satisfying $g^{-1}(z)<z(g^{-1})'(z)$ for each $z$. This is true provided $(g^{-1})''(z)$ is not zero in a neighbourhood of $0$. If it were, this would imply $g^{-1}$, and hence $g$, is linear in a neighbourhood of $0$. But $g(y)^n$ is concave by assumption, so this is impossible for $n\geq 2$. For similar technical reasons, it is false for $n=1$\footnote{For completeness, the argument is as follows: If $g(y)$ were a strictly increasing linear function, then the sum over the two points $x\in f^{-1}(s)$ of $1/|f'(x)|$ would be constant. But $1/|f'(x)|$ is a non-increasing function of $s$ when $x$ is the unique $x\in f^{-1}(s)$ greater than the critical point, and likewise for the unique less $ x\in f^{-1}(s)$ less than the critical point. Hence $|f'|$ is constant near the critical point, but if $f$ is $C^2$, this would be $0$, a contradiction.}.

Hence $nt(y)/t'(y)>y$ holds for all $y$ and $T$ has strictly increasing inverse, so $L\leq T^{-1}(\alpha)$. Thus $A_\alpha\subseteq\{x\in B:f(x)\leq\alpha T^{-1}(\alpha)\}$.

Set $\varepsilon_0=\alpha T^{-1}(\alpha)$. Then 
$$|J(\varepsilon_0)|=\left|\int_{\tilde{f}^{-1}(\varepsilon_0)}|\nabla \tilde{f}(x)|^{-1}\,d\mathcal{H}^{n-1}(x)\right|\leq\alpha^{-1}|\tilde{f}^{-1}(\varepsilon_0)|.$$
Since $\tilde{f}^{-1}(\varepsilon_0)$ is a sphere, we have $|\tilde{f}^{-1}(\varepsilon_0)|=n|S(\varepsilon_0)|^{(n-1)/n}$. Thus we have
$$|I(\lambda)|\leq 5n(\|a\|_\infty+\|a'\|_1)(|S(\varepsilon_0)|+(|\lambda|\alpha)^{-1}|S(\varepsilon_0)|^{(n-1)/n}).$$
Since $|S(\varepsilon_0)|=t(\alpha T^{-1}(\alpha))$, the bound becomes
$$5n(\|a\|_\infty+\|a'\|_1)t(\alpha T^{-1}(\alpha))$$
provided we choose $\alpha$ such that $t(\alpha T^{-1}(\alpha))=(|\lambda|\alpha)^{-1}t(\alpha T^{-1}(\alpha))^{(n-1)/n}$, that is, $(|\lambda|\alpha)^{-n}=t(\alpha T^{-1}(\alpha))$. By rearranging the definition of $T$, we have $t(xT^{-1}(x))=\omega_n(T^{-1}(x))^n$, so $\alpha$ must be chosen so that $\omega_n(T^{-1}(\alpha))^n=(|\lambda|\alpha)^{-n}$, that is,
$$\alpha=T\left(\frac{(|\lambda|\alpha)^{-1}}{\omega_n^{1/n}}\right)=\omega_n^{1/n}|\lambda|\alpha t^{-1}((|\lambda|\alpha)^{-n}),$$
equivalently,
$$t(\omega_n^{-1/n}|\lambda|^{-1})=(|\lambda|\alpha)^{-n}.$$
But $(|\lambda|\alpha)^{-n}=t(\alpha T^{-1}(\alpha))$, hence our bound is
$$|I(\lambda)|\leq 5n(\|a\|_\infty+\|a'\|_1)t(\omega_n^{-1/n}|\lambda|^{-1}).$$
This concludes the proof.

\textit{Acknowledgements.} The author is supported by a UK EPSRC scholarship at the Maxwell Institute Graduate School. The author would like to thank Prof. Jim Wright and Prof. Tony Carbery for helpful discussions and comments.

\bibliographystyle{abbrv}
\bibliography{Bibliography1}
John Green,\\ Maxwell Institute of Mathematical Sciences and the School of Mathematics,\\ University of Edinburgh,\\ JCMB, The King’s Buildings,\\ Peter Guthrie Tait Road,\\ Edinburgh, EH9 3FD,\\ Scotland\\ Email: \texttt{J.D.Green@sms.ed.ac.uk}
\end{document}